\documentclass[10pt, twocolumn]{article}
\usepackage{graphicx,amsmath,amssymb,url,enumerate,mathrsfs,epsfig,color}
\usepackage{tikz}
\usepackage{float}

\usepackage{amsmath}
\usepackage{amssymb}
\usepackage{yfonts}

\usepackage{tikz,pgfplots}

\usetikzlibrary{calc, patterns,arrows.meta, arrows, shapes.geometric}

\usepackage{graphicx,amsmath,amssymb,url,enumerate,mathrsfs,epsfig,color}

\usetikzlibrary{decorations.text}
\usetikzlibrary{decorations.markings}

\pgfplotsset{compat=1.8}

\usepackage{enumerate}
\usepackage{boxedminipage}
\usepackage{makeidx}
\usepackage{multicol}
\usepackage{xcolor}
\usepackage{mathtools}

\usetikzlibrary{calc, patterns,arrows, shapes.geometric}

\def\QED{\hskip0.1em\hfill\null\ \null\nobreak\hfill
\kern3pt\lower1.8pt\vbox{\hrule\hbox
{\vrule\kern1pt\vbox{\kern1.7pt \hbox{$\scriptstyle
QED$}\kern0.2pt}\kern1pt\vrule}\hrule}}

\newtheorem{thm}{Theorem}[section]

\newtheorem{corollary}{Corollary}[section]

\title{Hypercubes, $n$-groupoids, and mixtures }
\author{Marcelo Epstein\footnote{University of Calgary. Email: mepstein@ucalgary.ca}  }
\date{}
\begin{document}
\maketitle
\begin{abstract}
The theory of composite mixtures consisting of $n$ constituents is framed within the schema provided by the notion of $n$-groupoid. The point of departure is the analysis of $n$-dimensional hypercubes and their skeletons, to each of whose edges an element  (an arrow) of one of $n$ given material groupoids is assigned according to the coordinate class to which it belongs. In this way a $GL(3,{\mathbb R})$-weighted digraph is obtained. It is shown that if the double groupoid associated with each pair of constituents consists of commuting squares, the resulting $n$-groupoid is conservative. The core of this $n$-groupoid is transitive if, and only if, the mixture is materially uniform.
\end{abstract}


\section{Introduction}

The material response of a continuous medium is encapsulated in one or more {\it constitutive equations}. It is a remarkable fact that these equations endow the body manifold with an algebraic (and geometric) structure which can be best described as a {\it groupoid} \cite{ mackenzie,weinstein}, called the {\it material groupoid} associated with the given mechanical response.\footnote{This nomenclature appears first in \cite{jgp}, where the concept is extended to bodies with internal structure. Detailed treatments can be found in \cite{epsel, jofel, book}.} The basis for the construction of this groupoid lies in the notion of {\it material isomorphism}, which is roughly a linear map {\it or material transplant} between first-order neighbourhoods of two points rendering the material responses at them indistinguishable from each other. The arrows of this groupoid are precisely these material isomorphisms. The body is {\it materially uniform} (namely, made of the same material at all points) if, and only if, the material groupoid is transitive.

An extension of these ideas to certain kinds of {\it binary mixtures} was presented in \cite{symmetry, jmps}, neatly involving the more scantily explored applications of {\it double groupoids} \cite{natale, brown1}. It was found that even if the constituents of the mixture are perfectly uniform and defect-free, the mixture may fail to be so. Thus, a new kind of defectivity of {\it misalignment} was introduced. The present note attempts to promote a generalization of these ideas to mixtures with an arbitrary finite number $n$ of constituents by means of the concept of {\it $n$-groupoid} \cite{baez, zhu}. Physical applications of the theory of $n$-groupoids do not abound. It is hoped that this particular instance, with the concreteness associated to continuum mechanics, will contribute to at the very least trigger a useful physical imagery of the mathematical formulation, while also aiming at providing a basis for the definition of defects of misalignment in general mixtures.

\section{Hypercubes and skeletons}
\label{sec:hypercubes}

The solid $n$-dimensional unit hypercube $Q_n$ (or $n$-cube for short) is the geometric object defined as the Cartesian product
\begin{equation} 
Q_n=\underbrace{J \times ... \times J}_{\text{n~times}}= J^n \subset {\mathbb R}^n
\end{equation}
where $J=[0,1]$ is the oriented closed unit segment in $\mathbb R$.

For each integer $0 \le h <n$, an {\it $h$-face} $F_h$ is a subset of $Q_n$ obtained by fixing $n-h$ coordinates of ${\mathbb R}^n$ to either $0$ or $1$.  For example, a 4-face of $Q_7$ might be the set
\begin{equation}
F_4=\{0\}\times J \times\{1\}\times\{0\}\times J \times J \times J.
\end{equation}
Every $h$-face of $Q_n$ is congruent (by a Euclidean isometry) with the standard unit $h$-cube $Q_h$. For a fixed $h$, the number of $h$-faces is obtained as
\begin{equation}
N_h=2^{n-h} \,{n \choose h} =2^{n-h} \, \frac{n!}{h!\,(n-h)!}
\end{equation}

A $0$-face of $Q_n$ is called a {\it vertex}. The number of vertices is, therefore, $N_0=2^n$. A $1$-face is called an {\it edge}. The number of edges of $Q_n$ is $N_1=2^{n-1} n$. Thus, for example, the 3-dimensional cube has 12 edges. The edges of $Q_n$ can be partitioned into $n$ classes, according to the coordinate axis of ${\mathbb R}^n$ to which they  are parallel. If we denote the natural coordinates of ${\mathbb R}^n$ by $X_I$ (with $I=1,...,n$) we can loosely talk about the class of $X_I$-edges. For the ordinary unit cube in 3 dimensions (with coordinates $X,Y,Z$), we have 4 $X$-edges, 4 $Y$-edges, and 4 $Z$-edges.

An important relation between vertices is that of {\it adjacency}. Recalling that all vertices of $Q_n$ are points with coordinates in ${\mathbb R}^n$, each of which has a value 0 or 1, two vertices of $Q_n$ are said to be {\it adjacent with each other} if, and only if, their respective coordinates are identical except for a single entry. It follows from this definition that each vertex is adjacent with exactly $n$ other vertices. We remark that adjacency is not an equivalence relation since it is neither reflexive nor transitive. Since, by definition, an edge is obtained by fixing $n-1$ coordinates, edges join only adjacent vertices. The class to which an edge belongs is precisely determined by the coordinate $X_I$ in which the two adjacent vertices differ.

\begin{figure}[H]
\begin{center}
\begin{tikzpicture}[scale=0.8,
    line width=0.6pt,
    every node/.style={circle, draw, fill, minimum size=6pt, inner sep=0pt, font=\tiny\bfseries}]

    \pgfsetxvec{\pgfpoint{0.9cm}{0.0cm}}
    \pgfsetyvec{\pgfpoint{0.0cm}{0.9cm}}

    \foreach \point / \id / \angle in {
        (0,0)/1001/270,
        (0,5)/1011/90,
        (5,0)/1101/270,
        (5,5)/1111/90,
        (2,2)/0001/180,
        (2,7)/0011/90,
        (7,2)/0101/270,
        (7,7)/0111/90,
        (2.5,1.5)/1000/270,
        (2.5,3.5)/1010/180,
        (4.5,1.5)/1100/250,
        (4.5,3.5)/1110/120,
        (3.5,2.5)/0000/170,
        (3.5,4.5)/0010/180,
        (5.5,2.5)/0100/10,
        (5.5,4.5)/0110/0}
    {
        \node (\id) at \point [label=\angle:\id] {};
    }

    %

    \path 
    (0011) edge (0111) edge (1011) edge (0001)
    (0101) edge (0001) edge (1101) edge (0111)
    (1111) edge (1101) edge (0111) edge (1011)
   (0010) edge (0110) edge (1010) edge (0000)
   (0100) edge (0000) edge (1100) edge (0110)
    (1110) edge (1100) edge (0110) edge (1010)  ;

    \path[]
    (1001) edge (1101) edge (0001) edge (1011)
    (1000) edge (1100) edge (0000) edge (1010);
    
     \path[dotted]
    (0000) edge (0001)
    (0010) edge (0011)
    (0100) edge (0101)
    (0110) edge (0111)
    (1000) edge (1001)
    (1010) edge (1011)
    (1100) edge (1101)
    (1110) edge (1111);

\end{tikzpicture}
\end{center}
\caption{Tesseract obtained by joining corresponding vertices of two cubes by means of edges running along the fourth dimension (dotted lines)}
\label{fig:tesseract}
\end{figure}
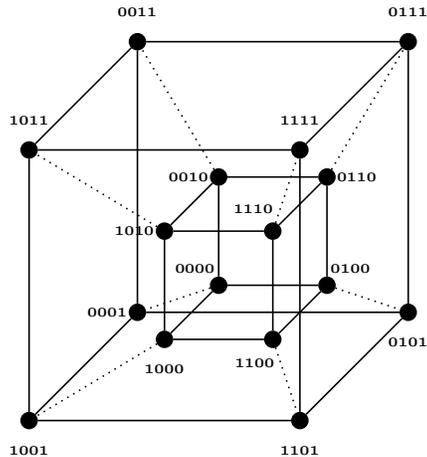

At the other extreme, for $h=n-1$, we obtain the so-called {\it facets} of the $n$-cube. The number of facets is $2 n$. There are exactly 2 facets perpendicular to each of the coordinate axes. We may think of $n$ classes of $X_I$-facets, each containing a pair of parallel facets. For example, the two $X_2$-facets of $Q_7$ are $ J  \times \{0\} \times J \times J \times J \times J \times J$ and $J  \times \{1\}\times J  \times J \times J \times J \times J$. Since each facet is congruent with $Q_{n-1}$, a creative and useful inductive way to think of $Q_n$ is as 2 parallel ($n-1$)-cubes of class $X_I$ in ${\mathbb R}^n$, separated by a unit distance along $X_I$ and joined by the $X_I$ unit segments between the corresponding vertices of each ($n-1$)-cube. In this way, we can think of the square as being generated by two parallel segments joined at their corresponding ends, and of a cube produced from two parallel squares joined at their corresponding vertices. To conceive of the 4-cube (or {\it tesseract}) we can imagine two ordinary cubes initially coincident, one of which is then subjected to a unit displacement into a putative fourth dimension, indicated with dotted lines in Figure \ref{fig:tesseract}.\footnote{Adapted from https://tex.stackexchange.com/questions/598573/} To complete the tesseract, 8 new edges are added running parallel to the fourth coordinate axis and connecting the corresponding vertices of the cubes.

\begin{figure}[]
\begin{center}
\begin{tikzpicture} [scale=1.0]
\tikzset{->/.style={decoration={
  markings,
  mark=at position .95 with {\arrow{stealth'}}},postaction={decorate}}}
 \draw[-stealth'] (0,0)--(4,0); 
  \draw[-stealth'] (0,0)--(0,4); 
   \draw[-stealth'] (0,0)--(-2,-2); 
   \node[below left] at (-2,-2) {$X^1$};
      \node[right] at (4,0) {$X^2$};
        \node[above] at (0,4) {$X^3$};
         \node[below left] at (-2,-2) {$X^1$};
\draw[thick,->] (0,0)--(3,0);
\draw[thick,->] (0,0)--(0,3);
\draw[thick,->] (3,0)--(3,3);
\draw[thick,->] (0,3)--(3,3);
\draw[thick,->] (0,3)--(-1.5,3-1.5);
\draw[thick,->] (3,3)--(3-1.5,3-1.5);
\draw[thick,->] (0,0)--(-1.5,-1.5);
\draw[thick,->] (3,0)--(3-1.5,-1.5);
\node at (0,0)[circle,fill,inner sep=2pt]{};
\node at (0,3)[circle,fill,inner sep=2pt]{};
\node at (3,0)[circle,fill,inner sep=2pt]{};
\node at (3,3)[circle,fill,inner sep=2pt]{};
\node at (-1.5,3-1.5)[circle,fill,inner sep=2pt]{};
\node at (3-1.5,3-1.5)[circle,fill,inner sep=2pt]{};
\node at (-1.5,-1.5)[circle,fill,inner sep=2pt]{};
\node at (3-1.5,-1.5)[circle,fill,inner sep=2pt]{};
\begin{scope}[xshift=-1.5cm,yshift=-1.5cm]
\draw[thick,->] (0,0)--(3,0);
\draw[thick,->] (0,0)--(0,3);
\draw[thick,->] (3,0)--(3,3);
\draw[thick,->] (0,3)--(3,3);
\end{scope}

\end{tikzpicture}
\end{center}
\caption{The oriented 3-skeleton $S_3$}
\label{fig:skeleton}
\end{figure}
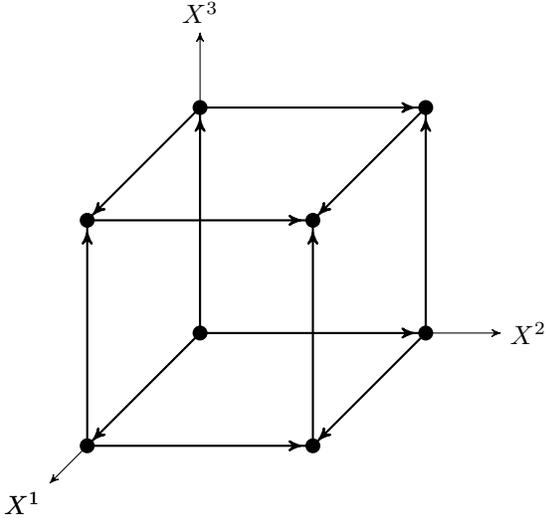

The construction just described can be made more precise by considering the {\it skeleton} $S_n$ of the $n$-cube $Q_n$. The skeleton $S_n$ is defined as the collection of vertices and edges of $Q_n$, occupying their positions in ${\mathbb R}^n$. This collection is also referred to as the {\it graph} of $Q_n$.\footnote{For a comprehensive treatment of this topic, see \cite{haray1, harary2}.} In our case, since we have considered an orientation of $J$, we obtain a {\it directed graph}, as shown in Figure \ref{fig:skeleton}. 

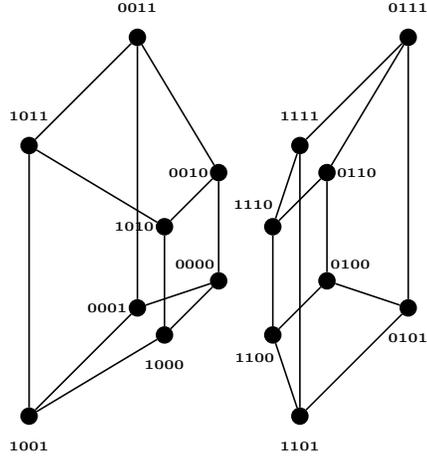
\begin{figure}[]
\begin{center}
\begin{tikzpicture}[scale=0.8,
    line width=0.6pt,
    every node/.style={circle, draw, fill, minimum size=6pt, inner sep=0pt, font=\tiny\bfseries}]

    \pgfsetxvec{\pgfpoint{0.9cm}{0.0cm}}
    \pgfsetyvec{\pgfpoint{0.0cm}{0.9cm}}

    \foreach \point / \id / \angle in {
        (0,0)/1001/270,
        (0,5)/1011/90,
        (5,0)/1101/270,
        (5,5)/1111/90,
        (2,2)/0001/180,
        (2,7)/0011/90,
        (7,2)/0101/270,
        (7,7)/0111/90,
        (2.5,1.5)/1000/270,
        (2.5,3.5)/1010/180,
        (4.5,1.5)/1100/250,
        (4.5,3.5)/1110/120,
        (3.5,2.5)/0000/170,
        (3.5,4.5)/0010/180,
        (5.5,2.5)/0100/10,
        (5.5,4.5)/0110/0}
    {
        \node (\id) at \point [label=\angle:\id] {};
    }

    \path 
   (0011) edge (1011) edge (0001)
    (0101)  edge (1101) edge (0111)
    (1111) edge (1101) edge (0111) 
   (0010)  edge (1010) edge (0000)
   (0100)  edge (1100) edge (0110)
    (1110) edge (1100) edge (0110)   ;

    \path[]
    (1001) edge (0001) edge (1011)
    (1000) edge (0000) edge (1010);
    
     \path[]
    (0000) edge (0001)
    (0010) edge (0011)
    (0100) edge (0101)
    (0110) edge (0111)
    (1000) edge (1001)
    (1010) edge (1011)
    (1100) edge (1101)
    (1110) edge (1111);

\end{tikzpicture}
\end{center}
\caption{The tesseract of Figure \ref{fig:tesseract} can also be obtained by joining corresponding vertices of the two cubes shown by means of edges (not shown) running along the second coordinate axis}
\label{fig:tesseract1}
\end{figure}

For a directed graph the notion of adjacency can be refined. Of two vertices adjacent with each other we declare that the vertex having an extra coordinate of value 1 is {\it adjacent from} the vertex having a corresponding 0 coordinate. Conversely, the latter vertex is said to be {\it adjacent to} the former. `Adjacency from' is not even a symmetric relation. The oriented skeleton $S_n$ consists of the vertices and the oriented edges (unit vectors) joining each vertex to those which are adjacent from it. The number of oriented edges issuing from or arriving at each vertex is not constant. Thus, the vertex at the origin is not adjacent from any other vertex. The inductive construction described above consistently delivers skeletons of all orders, without mixing solid ($n-1$)-cubes with one-dimensional edges. In terms of skeletons, therefore, it is easy to see that the $n$-skeleton generated by a pair of parallel ($n-1$)-skeletons is independent of the class of the two facets of departure. Indeed, the tesseract of Figure \ref{fig:tesseract} can be also obtained, as suggested in Figure \ref{fig:tesseract1}, by joining two ordinary cubes (in the space with coordinates $X^1, X^3, X^4$) by means of edges running along $X^2$.

\section{Coarse $n$-groupoids}

\subsection{Objective $n$-skeletons}
\label{sec:objective}

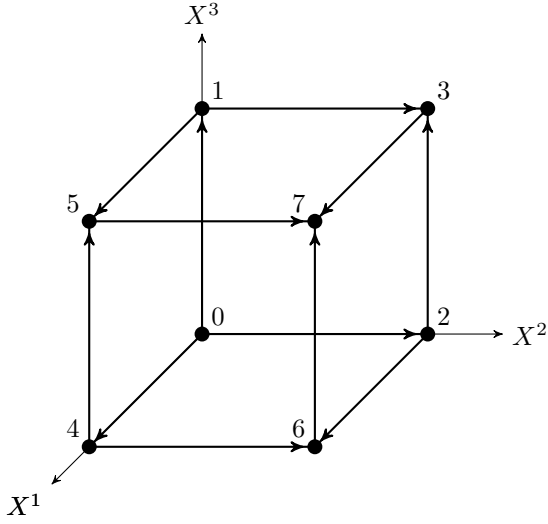
\begin{figure}[]
\begin{center}
\begin{tikzpicture} [scale=1.0]
\tikzset{->/.style={decoration={
  markings,
  mark=at position .95 with {\arrow{stealth'}}},postaction={decorate}}}
 \draw[-stealth'] (0,0)--(4,0); 
  \draw[-stealth'] (0,0)--(0,4); 
   \draw[-stealth'] (0,0)--(-2,-2); 
   \node[below left] at (-2,-2) {$X^1$};
      \node[right] at (4,0) {$X^2$};
        \node[above] at (0,4) {$X^3$};
         \node[below left] at (-2,-2) {$X^1$};
\draw[thick,->] (0,0)--(3,0);
\draw[thick,->] (0,0)--(0,3);
\draw[thick,->] (3,0)--(3,3);
\draw[thick,->] (0,3)--(3,3);
\draw[thick,->] (0,3)--(-1.5,3-1.5);
\draw[thick,->] (3,3)--(3-1.5,3-1.5);
\draw[thick,->] (0,0)--(-1.5,-1.5);
\draw[thick,->] (3,0)--(3-1.5,-1.5);
\node at (0,0)[circle,fill,inner sep=2pt]{};
\node at (0,3)[circle,fill,inner sep=2pt]{};
\node at (3,0)[circle,fill,inner sep=2pt]{};
\node at (3,3)[circle,fill,inner sep=2pt]{};
\node at (-1.5,3-1.5)[circle,fill,inner sep=2pt]{};
\node at (3-1.5,3-1.5)[circle,fill,inner sep=2pt]{};
\node at (-1.5,-1.5)[circle,fill,inner sep=2pt]{};
\node at (3-1.5,-1.5)[circle,fill,inner sep=2pt]{};
\begin{scope}[xshift=-1.5cm,yshift=-1.5cm]
\draw[thick,->] (0,0)--(3,0);
\draw[thick,->] (0,0)--(0,3);
\draw[thick,->] (3,0)--(3,3);
\draw[thick,->] (0,3)--(3,3);
\end{scope}
\node [above right] at (0,0) {$0$};
\node [above right] at (0,3) {$1$};
\node [above right] at (3,0) {$2$};
\node [above right] at (3,3) {$3$};
\node [above left] at (-1.5,-1.5) {$4$};
\node [above left] at (-1.5,3-1.5) {$5$};
\node [above left] at (3-1.5,-1.5) {$6$};
\node [above left] at (3-1.5,3-1.5) {$7$};

\end{tikzpicture}
\end{center}
\caption{Numbering the vertices of $S_3$ based on the binary system}
\label{fig:binary}
\end{figure}

Let $\mathcal B$ denote a set of objects and let $n>0$ be a natural number. We choose a {\it tuple} $\mathcal W$ of  $2^n$ points of $\mathcal B$ which we denote ${\mathcal W}=\{W_0, W_1,...,W_{2^n-1}\}$. Two or more points of $\mathcal W$ may not be necessarily distinct as objects of $\mathcal B$, but they are considered as different entries in $\mathcal W$. Since the unit $n$-cube has exactly $2^n$ vertices, it makes sense to define a bijection between $\mathcal W$ and the collection of vertices of $S_n$. For this purpose, we need to adopt a definite ordering of these vertices. One way to introduce a natural ordering of the vertices of ${\mathcal S}_n$ is to notice that the coordinates of each vertex of the unit $n$-cube consist of an ordered sequence of zeros and ones. But this is precisely the binary representation of the first $2^n$ non-negative integers! Denoting the ordered set of vertices of $S_n$ thus obtained by $\mathcal V=\{V_0, V_1,...,V_{2^n-1}\}$, as exemplified in Figure \ref{fig:binary}, we obtain the desired bijection $f:{\mathcal W} \to {\mathcal V}$ by $W_i \mapsto V_i$. Exploiting this natural bijection the various relations of adjacency can be pulled back to the tuple $\mathcal W$.

Introducing further structure, let us suppose that we have been given (on physical or mathematical grounds) $n$ simple groupoids ${\mathcal Z}_i \rightrightarrows {\mathcal B}\;(i=1,...,n)$ with common base $\mathcal B$. We denote by ${\mathcal P}^i_{XY}$ the set of arrows of $Z_i$ with source at $X \in {\mathcal B}$ and target at $Y \in {\mathcal B}$.\footnote{Unless the groupoid ${\mathcal Z}_i \rightrightarrows {\mathcal B}$ happens to be transitive, for some pairs $(X,Y)$ the set ${\mathcal P}^i_{XY}$ will be empty.} The inherited relation of `adjacency from' implies that each of the $2^n - 1$ entries $W_1,...,W_{2^n-1}$ of $\mathcal W$ is adjacent from at least one other element of $\mathcal W$ (including $W_0$). For each ordered pair $(W_a, W_b)$ of entries in $\mathcal W$ such that $W_b$ is adjacent from $W_a$ we choose an element of $P^I_{W_a W_b}$, where $I$ is determined by the coordinate $X_I$ in which $f(W_a)$ differs from $f(W_b)$. After all such pairs have been considered\footnote{If $P^I_{W_a W_b}$ is empty this construction is halted.} we obtain (with a blatant abuse of notation) an injective map $g:S_n \to {\mathcal W}\times {\mathcal P}_{{\mathcal W}\times{\mathcal W}}$. The image of this map (if not empty) is an {\it objective n-skeleton} ${\mathcal T}_n$ associated with $\mathcal W$. If each of the groupoids  ${\mathcal Z}_i \rightrightarrows {\mathcal B}\;(i=1,...,n)$ is transitive, at least one objective $n$-skeleton exists for each $2^n$-tuple $\mathcal W$. By having assigned to each objective oriented edge an arrow in a groupoid, we may regard the objective $n$-skeleton ${\mathcal T}_n$ as a kind of {\it weighted directed graph}, where the weights are arrows in various groupoids.

\subsection{Composition of objective n-skeletons}
\label{sec:composition}

Let $\tilde{\mathcal Z}$ denote the set of all objective $n$-skeletons arising from $n$ simple groupoids ${\mathcal Z}_i \rightrightarrows {\mathcal B}\;(i=1,...,n)$ with a common base $\mathcal B$ of objects. As demonstrated in Section \ref{sec:hypercubes}, the $n$-skeleton $S_n$ can be generated in $n$ different ways from any one of the $n$ pairs of opposite (parallel) facets. Because of the natural bijection between an objective $n$-skeleton ${\mathcal T}_n$ and the standard unit $n$-skeleton ${\mathcal S}_n$, these $n$ ways of conceiving an objective $n$-skeleton from pairs of parallel facets become available too for ${\mathcal T}_n$. This observation allows us to define an operation of composition  (denoted by $\odot_I$) in $\tilde{\mathcal Z}$ associated with each class of $X_I$-facets.

Let $\mathcal W$ be the $2^n$-tuple of vertices of an objective $n$-skeleton ${\mathcal T}_n$. The two $X_I$-facets of ${\mathcal T}_n$  will be denoted by ${\mathcal T}_{n-1}^{I,0}$ and ${\mathcal T}_{n-1}^{I,1}$. As noted in the description given in Section \ref{sec:hypercubes}, the corresponding vertices of these two facets are joined by the $2^{n-1}$ oriented $X_I$-edges of ${\mathcal T}_n$. Each of these edges is an arrow in the $I$-th groupoid ${\mathcal Z}_I \rightrightarrows {\mathcal B}$, with source and target projections $\alpha_I$ and $\beta_I$, respectively. Accordingly, we can define two projection maps
\begin{equation}
\tilde{\alpha}_I: \tilde{\mathcal Z} \to \tilde{\mathcal Z}_{\bar I}\;\;\;\;\;\;\;\;\tilde{\beta}_I: \tilde{\mathcal Z} \to \tilde{\mathcal Z}_{\bar I}
\end{equation}
where $\tilde{\mathcal Z}_{\bar I}$ denotes the set of all ($n-1$)-skeletons arising from the $n-1$ simple groupoids ${\mathcal Z}_i \rightrightarrows {\mathcal B}$ excluding the $I$-th groupoid. These two maps, called respectively the $I$-th source and target maps in $\tilde{\mathcal Z}$, are obtained by applying the source and target projections $\alpha_I$ and $\beta_I$ to each $X_I$-edge of an objective $n$-skeleton. Thus,
\begin{equation}
\tilde{\alpha}_I({\mathcal T}_n)= {\mathcal T}_{n-1}^{I,0}\;\;\;\;\;\;\;\;\tilde{\beta}_I({\mathcal T}_n)= {\mathcal T}_{n-1}^{I,1}
\end{equation}

Let ${\mathcal T}_n$ and ${\mathcal T}'_n$ be two $n$-skeletons in $\tilde{\mathcal Z}$ satisfying the condition
\begin{equation}
\tilde{\beta}_I ({\mathcal T}'_n) =\tilde{\alpha}_I ({\mathcal T}_n).
\end{equation}
Because of the way in which these projections were defined (in terms of the projections $\alpha_I$ and $\beta_I$), the corresponding $X_I$ edges of ${\mathcal T}_n$ and ${\mathcal T}'_n$ automatically satisfy the tip-to-tail condition as arrows of the $I$-th simple groupoid ${\mathcal Z}_I \rightrightarrows {\mathcal B}$ and they can be composed accordingly as such. We define the composition ${\mathcal T}_n  \odot_I  {\mathcal T}'_n$ as the objective $n$-skeleton such that
\begin{equation}
\tilde{\alpha}_I ({\mathcal T}_n  \odot_I  {\mathcal T}'_n)=\tilde{\alpha}_I( {\mathcal T}'_n),
\end{equation}
\begin{equation}
\tilde{\beta}_I ({\mathcal T}_n  \odot_I  {\mathcal T}'_n)=\tilde{\beta}_I( {\mathcal T}_n),
\end{equation}
and whose $X_I$-edges are obtained as the composition in the simple groupoid $Z_I$ of the corresponding $X_I$ edges of the factors.

Following along these lines (carefully defining the units, etcetera), we can endow the set $\tilde{\mathcal Z}$ with the structure of a simple groupoid over the set $\tilde{\mathcal Z}_{\hat I}$. Letting $I$ range from 1 to $n$ we obtain a collection of $n$ groupoid structures, one for each of the set of `objects' $\tilde{\mathcal Z}_{\hat I}$. The `arrows' of each of these simple groupoids consist of the set of $2^{n-1}$ oriented $X_I$-edges of an $n$-skeleton. Can these $n$ simple groupoid structures be unified under the umbrella of an $n$-groupoid? To obtain a positive answer to this question we need to show that the various compositions satisfy the compatibility conditions
\begin{equation}
 ({\mathcal T}_n  \odot_I  {\mathcal T}'_n) \odot_J  ({\mathcal T}''_n  \odot_I  {\mathcal T}'''_n) =  ({\mathcal T}_n  \odot_J  {\mathcal T}''_n) \odot_I  ({\mathcal T}'_n  \odot_J  {\mathcal T}'''_n)
 \end{equation}
for each distinct pair $I,J$ whenever the operations are possible. The proof is not difficult, but is omitted for now. We will refer to the set $\tilde{\mathcal Z}$ and its $n$ groupoid structures over $\tilde{\mathcal Z}_{\hat I}$ as an $n$-groupoid. The units and inverses of each $\tilde{\mathcal Z}_I$ are clearly defined by rendering the corresponding $X_I$-edges as units and inverses in the ordinary groupoid ${\mathcal Z}_I \rightrightarrows {\mathcal B}$.

Considering the recursive character of our definition of an $n$-groupoid, it should be clear that every 2-face of an $n$-groupoid is a double groupoid in the usual sense of the term \cite{ehresmann, natale, brown1}. To obtain the {\it core groupoid} of an $n$-groupoid we set all the edges, except those adjacent from the vertex $W_0$, of each objective $n$-skeleton to the unit of the corresponding groupoid.

\section{Material $n$-groupoids}

What we have obtained is the {\it coarse} $n$-groupoid of all possible $n$-skeletons constructed from the arrows of $n$ given simple groupoids. Imposing further restrictions on the groupoids of departure and on the selection of their arrows we can obtain smaller $n$-groupoids with a precise physical meaning, such as is the case in the theory of multi-composites and mixtures \cite{symmetry, jmps}.

Consider the special case in which $\mathcal B$ is a body manifold and each of the groupoids ${\mathcal Z}_i \rightrightarrows {\mathcal B}\;(i=1,...,n)$ is a subgroupoid of the 1-jet groupoid of $\mathcal B$, whose arrows are 1-jets of local diffeomorphisms between neighbourhoods of the points $X, Y \in {\mathcal B}$. In this case, the weights of the edges of our directed graphs become (in any given chart of $\mathcal B$) elements of the general linear group $GL(3, {\mathbb R})$. Because of these assumptions, the weights of two objective edges can be composed provided the target of an incoming edge is the same as the source of an outgoing edge from the same objective vertex. Using the language of graph theory, a {\it path} is a sequence of successively adjacent edges. In our weighted directed graph the weight of an edge traversed against its orientation is taken as the inverse of the original weight. The total {\it length} of the path is the composition of the weights of the traversed edges. A {\it circuit} is a path in which the first and last vertices coincide.

\begin{thm} The total weight of every circuit in ${\mathcal T}_n \in \tilde{\mathcal Z}$ vanishes if, and only if, every 2-face is commutative.
\end{thm}

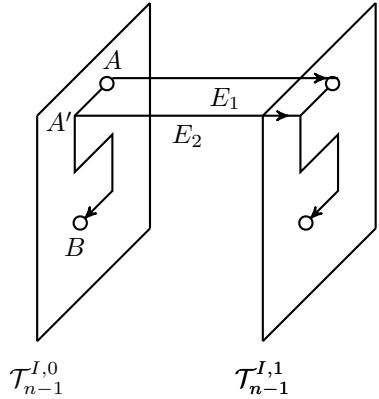
\begin{figure}[]
\begin{center}
\begin{tikzpicture} [scale=1.0]
\tikzset{->/.style={decoration={
  markings,
  mark=at position .95 with {\arrow{stealth'}}},postaction={decorate}}}

\draw[thick] (0,0)--(0,3);
\draw[thick] (0,3)--(-1.5,3-1.5);
\draw[thick] (0,0)--(-1.5,-1.5);
\draw[thick] (-1.5,-1.5)--(-1.5,1.5);
\draw[thick,o-o,->] (-0.5,2.)--(-1,1.5)--(-1.,0.75)--(-0.5,1.25)--(-0.5,0.5)--(-1,0);
\node at (-1.5,-2) {${\mathcal T}_{n-1}^{I,0}$};
\node[above] at (-0.5,2) {$A$};
\node[below] at (-1,0) {$B$};
\node[below left] at (-0.9,1.68) {$A'$};
\begin{scope} [xshift=3cm]
\draw[thick] (0,0)--(0,3);
\draw[thick] (0,3)--(-1.5,3-1.5);
\draw[thick] (0,0)--(-1.5,-1.5);
\draw[thick] (-1.5,-1.5)--(-1.5,1.5);
\draw[thick,o-o,->] (-0.5,2.)--(-1,1.5)--(-1.,0.75)--(-0.5,1.25)--(-0.5,0.5)--(-1,0);
\node at (-1.5,-2) {${\mathcal T}_{n-1}^{I,1}$};
\end{scope}
\node at (1.5,-2) {${\mathcal T}_{n-1}^{I,1}$};
\draw[thick, ->] (-0.5,2)--(2.5,2);
\draw[thick, ->] (-1,1.5)--(2.,1.5);
\node[below] at (1,2) {$E_1$};
\node[below] at (0.5,1.5) {$E_2$};

\end{tikzpicture}
\end{center}
\caption{Independence of weight for parallel traversing edges}
\label{fig:proof}
\end{figure}

We will refer to such a {\it conservative} $n$-groupoid as a {\it material n-groupoid}. The `only if' part of the theorem is trivial. The proof of the `if' part can be obtained by induction as follows. The theorem is trivially true for $n=2$. Assume it to be true for some value $n-1$. Given a circuit in ${\mathcal T}_n$,
\begin{enumerate}
\item If the circuit is completely contained in an ($n-1$)-face of ${\mathcal T}_n$ the theorem is obviously true.
\item If not, consider the partition of ${\mathcal T}_n$ into 2 $X_I$-facets. Since in a circuit the initial and terminal vertices coincide, any circuit must comprise an even number of $X_I$ edges, just like the number of bridge traversals between two shores connected by various bridges. Consider, therefore, a pair of successive traversals. If the connecting edge happens to be the same for both traversals, we conclude that their combined contribution to the total weight of the circuit vanishes since the same edge is traversed in opposite directions.
\item Otherwise, if the connecting edges are different, let $A$ and $B$ denote the corresponding vertices in the two $X_I$-facets joined by these two edges. We can join $A$ and $B$ by an identical finite path in each of the facets. Consider a single edge $AA'$ of these paths, as shown in Figure \ref{fig:proof}. Together with the two edges, say $E_1$ and $E_2$ joining the corresponding vertices, they constitute a 2-face of ${\mathcal T}_n$. By commutativity, therefore, the weight contributed by $E_1$ is identical to that contributed by $E_2$. Continuing along the path joining $A$ with $B$, we can invoke the reasoning adduced in the previous item (where the two edges coincided).
\end{enumerate}

\begin{corollary}  A mixture of $n$ uniform constituents is materially uniform if, and only if, the core of its material $n$-groupoid is transitive.
\end{corollary}

\section{Concluding remark}

If the mixture is not uniform we say that there are sources of misalignment, namely, a defectivity arising from an incompatibility between the constituents. We remark, however, that even if all the constituents are compatible with a fixed one, the mixture is not necessarily uniform. In other words, in the case of continuous symmetry groups, the core of the material $n$-groupoid may fail to be transitive, due to the extra degrees of freedom afforded by the symmetries. If this observation is correct, we can assert that a different kind of defectivity has been identified beyond that provided in \cite{jmps} for strictly binary mixtures.

\end{document}